\documentclass[12pt,a4paper,reqno]{amsart}
\usepackage{amssymb, amsbsy}

\newcommand{\No}{n.~}

\newtheorem{theor}{Theorem}

\newtheorem{state}[theor]{Proposition}
\newtheorem{lemma}[theor]{Lemma}
\theoremstyle{definition}

\theoremstyle{remark}
\newtheorem{cor}[theor]{Corollary}
\newtheorem{rem}{Remark}
\newtheorem{example}{Example}
\newtheorem{case}{Case}

\newcommand{\Id}{{\mathrm d}}

\newcommand{\pinner}{\mathbin{\mathchoice
   {\hbox{\vrule width0.6em depth0pt height0.4pt
   \vrule width0.4pt depth0pt height0.8ex}}
   {\hbox{\vrule width0.6em depth0pt height0.4pt
   \vrule width0.4pt depth0pt height0.8ex}}
   {\hbox{\kern0.14em
   \vrule width0.48em depth0pt height0.4pt
   \vrule width0.4pt depth0pt height0.6ex\kern0.14em}}
   {\hbox{\kern0.1em
   \vrule width0.39em depth0pt height0.4pt
   \vrule width0.4pt depth0pt height0.5ex\kern0.1em}}}}
\newcommand{\inner}{\pinner\,}

\DeclareFontFamily{OML}{cyr}{}
\DeclareFontShape{OML}{cyr}{m}{n}{
   <5> <6> <7> <8> <9> gen * wncyr
   <10> <10.95> <12> <14.4> <17.28> <20.74> <24.88> wncyr10
  }{}
\DeclareSymbolFont{rusletters}{OML}{cyr}{m}{n}
\DeclareSymbolFontAlphabet{\rusmath}{rusletters}
\DeclareMathSymbol\re{\rusmath}{rusletters}{"03}
\newcommand{\cEv}{\re}

\newcommand{\BBR}{{\mathbb{R}}}

\newcommand{\arctg}{\arctan}

\newcommand{\sym}{\mathop{\rm sym}\nolimits}

\newcommand{\const}{\mathop{\rm const}\nolimits}
\newcommand{\grad}{\mathop{\rm grad}\nolimits}

\newcommand{\bi}{{\boldsymbol{i}}}
\newcommand{\bs}{{\boldsymbol{s}}}

\newcommand{\bx}{{\boldsymbol{x}}}
\newcommand{\bE}{\mathbf{E}}

\newcommand{\cE}{\mathcal{E}}
\newcommand{\cEEL}{{\cE}_{\text{\textup{EL}}}}
\newcommand{\cEmS}{{\cE}_{{\min}\varSigma}}

\newcommand{\cL}{\mathcal{L}}
\newcommand{\LmS}{{L}_{{\min}\varSigma}}
\newcommand{\cLmS}{{\cL}_{{\min}\varSigma}}

\newcommand{\dd}{\partial}

\newcommand{\vph}{\varphi}

\newcommand{\even}{{\text{\textup{even}}}}
\newcommand{\odd}{{\text{\textup{odd}}}}

\newcommand{\ttbox}[1]{\mbox{\ttfamily#1}}
\newcommand{\by}[1]{\textit{{#1}}}
\newcommand{\jour}[1]{\textit{{#1}}}
\newcommand{\vol}[1]{\textbf{{#1}}}
\newcommand{\book}[1]{\textrm{{#1}}}

\title[On the minimal surface equation]{On the symmetry structure of
the minimal surface equation}
\date{October 26, 2004}

\author{A.~V.~Kiselev}
\thanks{A.\,K.\ was partially supported by the Lecce University
grant~\No\,650~CP/D}

\address{\textup{\textit{Permanent address}:}
153003 Russia, Ivanovo, Rabfakovskaya str.\ 34, Ivanovo State
Power University, Department of Higher Mathematics
\textup{and} Institute for Modelling and Computer Experiment.}
\curraddr{Department of Mathematics, Brock University,
500 Glenridge Ave.,  St.~Catharines, Ontario, Canada L2S~3A1.}
\email{Arthemy.Kiselev@brocku.ca}

\author{G.~Manno}
\address{Department of Mathematics, University of Lecce,
Via per Arnesano, 73100 Lecce (LE), Italy.}
\email{Gianni.Manno@kcl.ac.uk}

\thanks{Proc.\ conf.\ `\textit{Differential Geometry and Its
Applications}', Prague, Czech Republic, August~30~-- September~3,
2004.%
}

\subjclass[2000]{49Q05, 
                 53A10, 
                 58A20, 
                 70S10}  
\keywords{Minimal surface equation, symmetry, conservation laws}

\begin{document}
\rightline{ISPUmath-5/2004}

\begin{abstract}
An infinite sequence of commuting nonpolynomial contact symmetries of
the two\/-\/dimensional minimal surface equation is constructed.
Local and nonlocal conservation laws for $n$-dimensional minimal area
surface equation are obtained by using the Noether identity.
\end{abstract}
\maketitle

\subsection*{Introduction}
In this paper, we analyze the symmetry properties of the Euler equation
\begin{equation}\label{eqMinS}
\cEmS=\bigl\{\bE_u\bigl(\left[(1+(\grad u)^2\smash{)^{1/2}}\,\Id
\bx\right]\bigr)=0\bigr\}
\end{equation}
whose solutions describe the minimal area $n$-dimensional surfaces
$\varSigma=\{x^0=u(x^1$, $\ldots$, $x^n)\}\subset\BBR^{n+1}$.
We obtain variational conservation laws (either local or nonlocal) for
the equation $\cEmS$ in case $n\geq1$ is arbitrary and construct a
denumerable set of nonpolynomial contact symmetries of the minimal
surface equation if $n=2$ (note that the equation $\cEmS$ is the zero
mean curvature equation in this case).  We use the nonparametric
representation of the minimal surfaces $\varSigma$ such that at any
regular point within an open domain in the ambient manifold
$\BBR^{n+1}$ there is the fibre bundle $\pi\colon\BBR^{n+1}\to\BBR^n$
and any surface $\varSigma$ is the graph of a certain section in this
bundle. The minimal surface equation $\cEmS$ is a restriction upon
these sections. The latter equation itself is a submanifold
$\cEmS\subset J^2(\pi)$ of the second jet space of the fibre bundle
$\pi$, see~\cite{ClassSymEng} for notation,
definitions, and details. From now
on, we assume that the symmetries, conservation laws, and other
structures are restricted onto the infinite prolongation
$\cEmS^\infty\subset J^\infty(\pi)$; the superscipt $\infty$ will be
omitted by default.

The analytic and topological properties of the minimal surfaces
(\textit{i.e.}, either the minimal area surfaces of the zero mean
curvature surfaces) are considered in the vastest
literature~(\cite{Osserman}).
Recently, classical results and concepts on this topic were enlarged by
their applications in mathematical physics (\textit{e.g.}, the harmonic
maps theory and the $\sigma$-models, see~\cite{FordyWood}).
Nevertheless, the symmetry structure of the minimal surface equation
itself remains rather unveiled. In the paper~\cite{Bila}, seven
polynomial symmetries of the two\/-\/dimensional equation $\cEmS$
(namely, the shift, two translations, three rotations, and the
dilatation) were obtained~--- and a denumerable set of nonpolynomial
contact symmetries~(\cite{ISPU}) was missed at all.

The aim of this
paper is to fill in the apparent gap in description
of the symmetry algebra
for the two\/-\/dimensional minimal surface equation.
Also, we reconstruct several conservation laws for the $n$-dimensional
equation~$\cEmS$.
The knowledge of local and nonlocal conservation laws for the minimal
surface equation is expected to contribute to further studies of the
recursion operators for (not necessarily local part of) the symmetry
algebra of this equation.

\section{Conservation laws of the $n$-dimensional equation}
\noindent%
First, we recall known symmetries of the minimal area surface equation.

\begin{state}[\textup{\cite{Bila}}]\label{BilaSym}
The generators of the Lie algebra of point symmetries for the
two\/-\/dimensional minimal surface equation
\begin{equation}\label{eqMinS2}
\cEmS=\{(1+u_y^2)u_{xx}-2u_xu_yu_{xy}+(1+u_x^2)u_{yy}=0\}
\end{equation}
are
%
\begin{align*}
\vph_1&=1 & &\text{\textup{(}the shift\textup{),}}\\
\vph_2^i&=u_{x^i} & &\text{\textup{(}translations\textup{),}}\\
\begin{aligned}\vph_3^{12}&{}\\ \bar\vph_3^i&{}\end{aligned}
&
\begin{aligned}{}&\!=yu_{x}-xu_{y}\\ {}&\!=x^i+uu_{x^i}\end{aligned}
& &\text{\textup{(}rotations\textup{),}}\\
\vph_4&=u-\langle \grad u, \bx\rangle
   & &\text{\textup{(}the dilatation\textup{).}}
\end{align*}
\end{state}

\begin{rem}\label{RemSymNDim}
The sections $\vph_1$, $\vph_2^i$,
$\vph_3^{ij}=x^j u_{x^i}- x^i u_{x^j}$,
$\bar\vph_3^i$, 
and $\vph_4=u-\langle\grad u$, $\bx\rangle$ are symmetries of
the $n$-dimensional minimal area surface Euler equation
$\cEmS=\{\bE_u(\cLmS)=0\}$ if $n\geq1$ is arbitrary and $1\leq i$,
$j\leq n$; here the Lagrangian is $\cLmS=[\LmS\,\Id\bx]$ and
$\LmS=\sqrt{1+(\grad u)^2}$.

Also, we note that the sections
$\vph_1$, $\vph^i_2$, $\vph_3^{ij}$, and
$\bar\vph_3^i$ generate
the $\tfrac{1}{2}(n+1)(n+2)$-di\-men\-si\-o\-nal
algebra of the Killing vector fields
on the Euclidean space~$\BBR^{n+1}\supset\varSigma$.
\end{rem}

Now we assign continuity equations to the symmetries that preserve the
Lagrangian.  

\begin{rem}
The resulting conserved currents are
either purely local or involve nonlocalities. This is owing to three
possible ways for a Lagrangian~$\cL$
to be conserved w.r.t.\ the symmetries of
the corresponding Euler equation~$\cEEL$:
\end{rem}

\begin{itemize}
\item
If the functional $\cL$ is conserved on the whole jet space
$J^\infty(\pi)\supset\cEEL^\infty$ along a symmetry $\vph\in\sym\cEEL$,
then $\vph$ is a \emph{Noether symmetry} of the Lagrangian~$\cL$ and the
local conservation law is assigned to $\vph$ by the first Noether
theorem~(\cite{Mystique, ClassSymEng}). Then we use the following lemmas.
\end{itemize}

\begin{lemma}[\textup{Noether's identity, \cite{Ibragimov, Vin2Cal}}]
Let $\cEEL=\{F\equiv\bE_u(\cL)=0\}$ be the Euler equation assigned to a
Lagrangian $\cL=[L\,\Id\bx]$
and let $\vph$ be a Noether symmetry\textup{:}
$\cEv_\vph(\cL)=\Id_h(\mu)$,
where $\Id_h=\sum_{i=1}^n \Id x^i\otimes D_i$
is the standard horizontal differential.
Then $\vph$ is the generating section of
the conservation law $[\omega]=[\mu-\nu-\lambda]$\textup{:}
$\Id_h(\omega)=\square(F)\,\Id\bx$ such that
$\vph=\square^*(1)$\textup{;}
here $\cEv_\vph(\cL)=\langle\vph,F\rangle+\Id_h(\nu)$ and
$\langle1,\square(F)\rangle=\langle\vph,F\rangle+\Id_h(\lambda)$.
\end{lemma}

\begin{proof}
We have
\[
\Id_h(\mu)=\cEv_\vph(\cL)=\langle1,\ell_L^*(\vph)\rangle=
\langle\ell_L^*(1),\vph\rangle+\Id_h(\nu),
\]
where $\ell_L^*(1)=\bE_u(\cL)$ by definition.
Therefore, $\Id_h(\mu-\nu)=\langle1,\square(F)\rangle+\Id_h(\lambda)$
such that $\vph=\square^*(1)$.
\end{proof}

\begin{lemma}[\textup{\cite{Ibragimov}}]
In local coordinates, the Noether identity is
\begin{equation}\label{NoetherCoords}
\cEv_\vph=\vph\cdot\bE_u+\sum\nolimits_{i=1}^n D_i\circ Q_{\vph,i},
\end{equation}
where
\[
Q_{\vph,i}=\sum_\tau\sum_{\rho+\eta=\tau}(-1)^\eta D_\rho(\vph)\cdot
D_\eta\circ\dd/\dd u_{\tau+1_i}.
\]
\end{lemma}

\begin{proof}
Fix a multi\/-\/index $\sigma$; then the coefficient of $\dd/\dd u_\sigma$
in the r.h.s.\ of Eq.~\eqref{NoetherCoords} is
\begin{multline*}
\vph\cdot(-1)^\sigma\,D_\sigma+\sum_{i=1}^n D_i\circ\Bigl[
\sum_{\tau+1_i=\sigma}\sum_{\rho+\eta=\tau}(-1)^\eta\,D_\rho(\vph)\cdot
D_\eta\Bigr] = {}\\
{}=\vph\cdot(-1)^\sigma\,D_\sigma+
\sum_{\substack{\rho'+\eta=\sigma\\ |\rho'|>0}}(-1)^\eta\,
D_{\rho'}(\vph)\cdot D_\eta -
\sum_{\substack{\rho+\eta'=\sigma\\ |\eta'|>0}}(-1)^{\eta'}\,
D_{\rho}(\vph)\cdot D_{\eta'}={}\\
{}=\vph\cdot(-1)^\sigma\,D_\sigma+D_\sigma(\vph)-(-1)^\sigma\,\vph\cdot
D_\sigma=D_\sigma(\vph).
\end{multline*}
This completes the proof.
\end{proof}

\begin{itemize}
\item
If $\cL$ is conserved on the equation $\cEEL$ only (\textit{i.e.},
$[\cEv_\vph(\cL)]=\const\cdot\cL$),
then $\vph$ is a \emph{variational symmetry} of
$\cEEL$ and the conservation law reconstructed by $\vph$ involves
nonlocalities.
\end{itemize}

\noindent%
Indeed, we have
\[
\Id_h(\mu)+\const\cdot\cL=\cEv_\vph(\cL)=
\langle\bE_u(\cL),\vph\rangle+\Id_h(\nu)=
\langle1,\square(F)\rangle+\Id_h(\nu+\lambda),
\]
where $\vph=\square^*(1)$.
Therefore, if a set of nonlocalities $\bs$ trivializes the Lagrangian:
$\tilde\Id_h(\eta)=\cL$,
here we put $\eta=\sum_i s^i\cdot(\dd/\dd x^i\inner\Id\bx)$,
then we obtain a nonlocal conservation law $\Omega=\mu+\eta-\nu-\lambda$
whose generating section is~$\vph$.
Generally, the compatibility conditions $\bs_{x^ix^j}=\bs_{x^jx^i}$
determine an infinite set of new variables.

\begin{itemize}
\item
Finally, if the Lagrangian $\cL$ is not conserved neither on the jet
space $J^\infty(\pi)$ nor on the equation $\cEEL$, then there is
no conservation law for the symmetry~$\vph$.
\end{itemize}

We obtain the following assertion.

\begin{state}\indent\par
\begin{enumerate}
\item
Under the notation of Proposition~\textup{\ref{BilaSym}} and
Remark~\textup{\ref{RemSymNDim}},
the point symmetries $\vph_1$, $\vph_2$, and
$\vph_3^{ij}$ are Noether's symmetries of the Lagrangian $\cLmS$.
Indeed, the equality $[\cEv_{\vph_k}(\cLmS)]=0$ holds on the jet space
$J^\infty(\pi)$ for $1\leq k\leq 3$.
The conserved currents assigned to these symmetries are such that the
corresponding continuity equations are, respectively,
\begin{gather*}
\sum\nolimits_{i=1}^n \bar D_i(\dd L/\dd u_{x^i})=0,\\
\bar D_i\left(u_{x^i}\cdot\dd L/\dd u_{x^i}-L\right) +
\sum\nolimits_{j\not=i}\bar D_j\left(u_{x^i}\cdot\dd L/\dd u_{x^j}\right)=0,\\
%
\bar D_i(x^jL-\vph_3^{ij}u_{x^i}/L) -
\bar D_j(x^iL+\vph_3^{ij}u_{x^j}/L)
{}-\sum\nolimits_{k\not=i,j} \bar D_k(\vph_3^{ij}\cdot u_{x^k}/L)
=0.
\end{gather*}
\item
The dilatation $\vph_4$ preserves the Lagrangian $\cLmS$ on the
equation $\cEmS$\textup{:}
\[
[\cEv_{\vph_4}(\cLmS)]=n\,\cLmS,
\]
and hence
\[
\bE_u(\cEv_{\vph_4}(\cLmS))=n\,\bE_u(\cLmS)=0\quad\text{on\quad
$\cEmS=\{\bE_u(\cLmS)=0\}$.}
\]
Let $s^i$, $1\leq i\leq n$, be nonlocal variables such that their
derivatives be $s^i_{x^i}=\LmS$. Then the nonlocal
continuity equation
\[
\sum\nolimits_i \tilde D_i \Bigl( s^i - x^i\,L -
\frac{uu_i\mathstrut}{\sqrt{1+\sum\nolimits_j u_j^2}\mathstrut} +
\frac{u_i\cdot\sum\nolimits_j x^j u_j \mathstrut}%
{\sqrt{1+\sum\nolimits_{k\mathstrut} u_k^2}\mathstrut}
\Bigr) = 0
\]
is assigned to the symmetry~$\vph_4$.
\end{enumerate}
\end{state}

\noindent%
The proof is straightforward.

\section{Contact symmetries of the two\/-\/dimensional equation}
\noindent%
From now on we assume $n=2$. We claim that the symmetries $\vph_1$,
$\ldots$, $\vph_4$ do not exhaust the set of generators of the
contact symmetry Lie algebra for the minimal surface equation~$\cEmS$,
see Eq.~\eqref{eqMinS2}. Indeed, we have

\begin{state}
The Lie algebra of contact symmetries for the minimal surface equation
$\cEmS$ is generated by the solutions $\vph(u_x$\textup{,} $u_y)$ of
the equation
\begin{equation}\label{eqGianni}
(1+u_x^2)\,\frac{\dd^2\vph}{\dd(u_x)^2} +
2u_xu_y\,\frac{\dd^2\vph}{\dd u_x\dd u_y} +
(1+u_y^2)\,\frac{\dd^2\vph}{\dd(u_y)^2} = 0
\end{equation}
\textup{(}\textit{e.g.}, the shift $\vph_1=1$ and the translations
$\vph_2^i=u_{x^i}$ are solutions to Eq.~\eqref{eqGianni}\textup{)}
and the sections $\vph_3^{12}$, $\bar\vph_3^i$, and $\vph_4$\textup{;}
here $i=1$\textup{,}\,$2$.
\end{state}

\noindent%
The proof is elementary by using the analytic transformations
software~(\cite{Jet}) and therefore omitted. We also note that
Eq.~\eqref{eqGianni} and its solutions except the point symmetries
$\vph_1$ and $\vph_2^i$ were missed in~\cite{Bila}.

Now we construct a denumerable set (\cite{ISPU}) contact symmetries
of the minimal surface equation~\eqref{eqMinS2}.
Consider Eq.~\eqref{eqGianni} and suppose that $\vph(u_x$, $u_y)$ is
polynomial in $u_y$ (of course, all reasonings are preserved by the
symmetry transformation $x\leftrightarrow y$).
Then we obtain two distinct cases: the degree $K$ of the polynomial at
hand can be either even ($K=2k$) or odd ($K=2k+1$). In what follows, we
treat these cases separately.

\begin{case}[\textup{$K=2k$}]
We assume $\vph^\even_k=\sum_{\ell=0}^k f_\ell(u_x)\,u_y^{2\ell}$.
Then from Eq.~\eqref{eqGianni} we obtain the following chain of
equations:
\begin{enumerate}
\item
The homogeneous equation
\begin{equation}\label{eqEvenTop}
(1+u_x^2)\,f_k'' + 4ku_x\,f_k' + 2k(2k-1)\,f_k = 0
\end{equation}
at the highest power $u_y^{2k}$.
\item
The nonhomogeneous equations
\begin{equation}\label{eqEvenIntermed}
(1+u_x^2)\,f_\ell'' + 4\ell u_x\,f_\ell' + 2\ell(2\ell-1)\,f_\ell =
-(2\ell+2)(2\ell+1)\,f_{\ell+1}
\end{equation}
at the intermediate powers $u_y^{2\ell}$, $1\leq \ell<k$;
here the homogeneous components are Eq.~\eqref{eqEvenTop} with the
index shift.
\item
The terminal relation at the zero power~$u_y^0$:
\begin{equation}\label{eqEvenBottom}
(1+u_x^2)\,f_0'' = -2\,f_1.
\end{equation}
\end{enumerate}
The solutions of system~(\ref{eqEvenTop}--\ref{eqEvenBottom}) are the
following.
\begin{enumerate}
\item
The solution of homogeneous equation~\eqref{eqEvenTop} is
$f_k(u_x)=C_1\cdot(\bi+u_x)^{-2k+1} + C_2\cdot(-\bi+u_x)^{-2k+1}$,
where the constants $C_1$, $C_2$ are arbitrary.
The equivalent basis of real solutions to Eq.~\eqref{eqEvenTop} is
\begin{align*}
f_{k}^1&={(1+u_x^2)^{-2k+1}}\cdot
  \sum_{\ell=0}^{k-1} (-1)^\ell\binom{2k-1}{2\ell} u_x^{2k-2\ell-1},\\
f_{k}^2&={(1+u_x^2)^{-2k+1}}\cdot
  \sum_{\ell=1}^{k} \binom{2k-1}{2\ell-1} u_x^{2k-2\ell}.
\end{align*}
\item
The solutions to the intermediate nonhomogeneous
equations~\eqref{eqEvenIntermed} are
\begin{multline*}
f_\ell(u_x)={}\\
(2\ell+2)(2\ell+1)\cdot\Bigl\{
  f_\ell^1(u_x)\cdot\int\frac{f_{\ell+1}\,f_\ell^2\,\Id u_x}{f_\ell^1\,{(f_\ell^2)}' -
  {(f_\ell^1)}'\,f_\ell^2} -
  f_\ell^2(u_x)\cdot\int\frac{f_{\ell+1}\,f_\ell^1\,\Id u_x}{f_\ell^1\,{(f_\ell^2)}' -
  {(f_\ell^1)}'\,f_\ell^2}
\Bigr\}.
\end{multline*}
All integration constants can be omitted since they originate from the
homogeneous equations that provide the symmetries polynomial in $u_y$
whose degrees are less than $K=2k$.
\item
The quadrature for terminal equation~\eqref{eqEvenBottom} is
\[
f_0(u_x)=-2\int^{u_x}\Id\xi\int^\xi\frac{f_1(\eta)\,\Id\eta}{1+\eta^2}
+\alpha u_x+\beta,\qquad\alpha,\beta\in\BBR.
\]
The coefficients of the integration constants $\alpha$, $\beta$ provide
the translation $u_x$ and the shift $1\in\sym\cEmS$, respectively.
\end{enumerate}
\end{case}

\begin{case}[\textup{$K=2k+1$}]
Now we have $\vph^\odd_k=\sum_{\ell=0}^k g_\ell(u_x)\,u_y^{2\ell+1}$.
From Eq.~\eqref{eqGianni} we obtain the chain of equations:
\begin{enumerate}
\item
The homogeneous equation
\begin{equation}\label{eqOddTop}
(1+u_x^2)\,g_k'' + 2(2k+1)ku_x\,g_k' + 2k(2k+1)\,g_k = 0
\end{equation}
at the highest power $u_y^{2k+1}$.
\item
The nonhomogeneous equations
\begin{equation}\label{eqOddIntermed}
(1+u_x^2)\,g_\ell'' + 2(2\ell+1)u_x\,g_\ell' + 2\ell(2\ell+1)\,g_\ell =
 -(2\ell+3)(2\ell+2)\,g_{\ell+1}
\end{equation}
at the intermediate powers $u_y^{2\ell+1}$, $1\leq\ell<k$;
again, the homogeneous components are Eq.~\eqref{eqOddTop} with the
index shift.
\item
The terminal relation at~$u_y$:
\begin{equation}\label{eqOddBottom}
(1+u_x^2)\,g_0'' + 2u_x\,g_0' = -6\,g_1.
\end{equation}
\end{enumerate}
The solutions of system~(\ref{eqOddTop}--\ref{eqOddBottom}) are the
following.
\begin{enumerate}
\item
The solution of homogeneous equation~\eqref{eqOddTop} is
\[
g_k(u_x)=C_3\cdot(1+u_x^2)^{-k}\cdot
\left(\frac{\bi+u_x}{-\bi+u_x}\right)^k +
  C_4\cdot(1+u_x^2)^{-k}\cdot
\left(\frac{-\bi+u_x}{\bi+u_x}\right)^k,
\]
where $C_3$, $C_4$ are arbitrary constants.
The equivalent basis of real solutions to Eq.~\eqref{eqOddTop} is
\begin{align*}
g_{k}^1&={(1+u_x^2)^{-2k}}\cdot
  \sum_{\ell=0}^{k-1} (-1)^\ell\binom{2k}{2\ell} u_x^{2k-2\ell},\\
g_{k}^2&={(1+u_x^2)^{-2k}}\cdot
  \sum_{\ell=1}^{k} \binom{2k}{2\ell-1} u_x^{2k-2\ell+1}.
\end{align*}
\item
The solutions to the intermediate nonhomogeneous
equations~\eqref{eqOddIntermed} are
\begin{multline*}
g_\ell(u_x)={}\\
(2\ell+3)(2\ell+2)\cdot\Bigl\{
  g_\ell^1(u_x)\cdot\int\frac{g_{\ell+1}\,g_\ell^2\,\Id u_x}{g_\ell^1\,{(g_\ell^2)}' -
  {(g_\ell^1)}'\,g_\ell^2} -
  g_\ell^2(u_x)\cdot\int\frac{g_{\ell+1}\,g_\ell^1\,\Id u_x}{g_\ell^1\,{(g_\ell^2)}' -
  {(g_\ell^1)}'\,g_\ell^2}
\Bigr\}.
\end{multline*}
Again,
the integration constants can be omitted since they originate from the
homogeneous equations that provide the symmetries
which are polynomial in $u_y$ and whose degrees are less than $K=2k+1$.
\item
The quadrature for terminal equation~\eqref{eqOddBottom} is
\[
g_0(u_x)=-6\int^{u_x}\frac{\Id\xi}{1+\xi^2}\int^\xi g_1(\eta)\,\Id\eta
+ \gamma + \delta\arctan u_x,\qquad\gamma,\delta\in\BBR.
\]
The coefficients of the constants $\gamma$, $\delta$ are the
translation $u_y$ and the first nonpolynomial symmetry
$\vph_5=u_y\arctan u_x$, respectively.
\end{enumerate}
\end{case}

The above reasonings provide two sequences of the symmetries
$\vph_k^{\even,\,\odd}\in\sym\cEmS$
that are polynomial in $u_y$ and nonpolynomial
(except three starting terms) in~$u_x$.

\begin{example}\label{ExampleStart}
The initial terms of these sequences are
\begin{align*}
\vph_1&=1, \quad \vph_2^1=u_x, &
\vph_2^2&=u_y, \quad \vph_5=u_y\arctg u_x,\\
\vph_6&=\frac{u_xu_y^2}{1+u_x^2}+\arctg u_x, &
  \vph_{7}&=\frac{u_y^2}{1+u_x^2}-u_x\arctg u_x,\\
\vph_{8}&=\frac{u_xu_y^3}{(1+u_x^2)^2}+\tfrac{3}{2}\frac{u_xu_y}{1+u_x^2}, &
  \vph_{9}&=\frac{u_x^2-1}{(1+u_x^2)^2}\cdot
  u_y^3-\frac{3u_y}{1+u_x^2}.
\end{align*}
\end{example}

\begin{rem}
Analogous reasonings can be applied to the formal power series
solutions
\[
\vph_k^\even=\sum_{\ell\leq-k} f_\ell(u_x)\,u_y^{2\ell},\qquad
\vph_k^\odd=\sum_{\ell\leq-k} g_\ell(u_x)\,u_y^{2\ell+1},
\]
where $k>0$. The symmetries $\vph_k^\even$ are therefore polynomial
in~$u_y^{-1}$.
\end{rem}

\begin{rem}
There are two different types of recursion operators that act on the
contact symmetries $\vph_k(u_x$, $u_y)$. Let $k\geq0$, $i=1$,\,$2$, and
suppose $\vph_{k,i}$ is the symmetry such that $f_k^i$ (resp., $g_k^i$)
is the coefficient of the highest power of $u_y$. Then we have
\begin{itemize}
\item
the operator $\Delta\colon\vph_{k,1}\rightleftarrows\vph_{k,2}$ that
swaps the elements of the basis if the highest power that depends on
$k$ is fixed;
\item
the operator $\nabla\colon\vph_{k,i}\mapsto\vph_{k+1,i}$ that
proliferates the symmetry heads along the sequence by the rules
\begin{align*}
\nabla^\even\colon f_k^1&\mapsto\frac{1}{(\bi+u_x)^2}\,f_k^1=f_{k+1}^1,
&
\nabla^\even\colon
f_k^2&\mapsto\frac{1}{(-\bi+u_x)^2}\,f_k^2=f_{k+1}^2,
\\
\nabla^\odd\colon g_k^1&\mapsto\frac{-1}{(\bi+u_x)^2}\,g_k^1=g_{k+1}^1,
&
\nabla^\odd\colon
g_k^2&\mapsto\frac{1}{(-\bi+u_x)^2}\,g_k^2=g_{k+1}^2.
\end{align*}
\end{itemize}
\end{rem}

\begin{lemma}
Assume that $\vph'(u_x, u_y)$ and $\vph''(u_x, u_y)$ are the generating
sections of evolutionary vector fields $\cEv_{\vph'}$
and~$\cEv_{\vph''}$. Then their Jacobi bracket $\{\vph'$\textup{,}
$\vph''\}$ is always trivial.
\end{lemma}
\begin{proof}
In local coordinates, the Jacobi bracket $\{\vph'$\textup{,}
$\vph''\}=\cEv_{\vph'}(\vph'')-\cEv_{\vph''}(\vph')$ is
\begin{multline*}
\{\vph',\vph''\}=
 \frac{\dd\vph'}{\dd u_x}u_{xx}\cdot\frac{\dd\vph''}{\dd u_x}+
 \frac{\dd\vph'}{\dd u_y}u_{xy}\cdot\frac{\dd\vph''}{\dd u_x}+
 \frac{\dd\vph'}{\dd u_x}u_{xy}\cdot\frac{\dd\vph''}{\dd u_y}+
 \frac{\dd\vph'}{\dd u_y}u_{yy}\cdot\frac{\dd\vph''}{\dd u_y} -
\text{v.\,v.} = 0.
\end{multline*}
\end{proof}

\begin{cor}
The contact symmetries $\vph_k(u_x$, $u_y)$ of the minimal surface
equation~$\cEmS$ commute.
\end{cor}

\begin{state}
Suppose that a minimal surface $\varSigma$ is invariant w.r.t.\ a
contact symmetry $\vph(u_x$, $u_y)$. Then $\varSigma$ is a plane.
\end{state}
\begin{proof}
Consider the constraint $\vph=0$. By the implicit function theorem, we
have $u_y=\phi(u_x)$ almost everywhere. Therefore
$u_{yy}=\left(\phi'(u_x)\right)^2\cdot u_{xx}$ and from
Eq.~\eqref{eqMinS2} we obtain the equation
\[
\Bigl(
(1+u_x^2)\cdot\left(\phi'(u_x)\right)^2 -
2u_x\,\phi(u_x)\cdot \phi'(u_x) + (1+\phi^2(u_x))\Bigr) \cdot u_{xx} = 0.
\]
Hence either $u_{xx}=0$ and we have $u_x=\const$, $u_y=\phi(u_x)=\const$,
or $u_x$ is subject to the algebraic equation such that its solutions
are $u_x=\const$ and therefore $u_y=\phi(u_x)=\const$ again.
\end{proof}

We conjecture that none of the contact symmetries $\vph_k$, $k\geq5$
(see Example~\ref{ExampleStart}), are Noether.

\subsection*{Acknowledgements}
The authors thank I.\,S.\,Krasil'shchik, A.\,M.\,Verbovetsky,
and R.\,Vitolo for useful discussions.
A.\,K.\ is grateful to V.\,Rosenhaus for stimulating
remarks and to the University of Lecce for hospitality. The research of
A.\,K.\ was supported by the University of Lecce grant \No650~CP/D.

\end{document}